\input amstex.tex
\documentstyle{amsppt}
\magnification=1200
\loadbold

\def\Co#1{{\Cal O}_{#1}}
\def\Fi#1{\Phi_{|#1|}}
\def\fei#1{\phi_{#1}}

\def\lrw{\longrightarrow}

\def\Bbbp1{{\Bbb P}^1}
\def\simlin{\sim_{\text{\rm lin}}}
\def\simnum{\sim_{\text{\rm num}}}

\def\Div{\text{\rm Div}}

\def\dim{\text{\rm dim}}
\def\roundup#1{\ulcorner{#1}\urcorner}
\TagsOnLeft
%%%%%%%%%%%%%%%%%%%%%%%%%%%
\topmatter
\title
On pluricanonical maps for threefolds of general type, II
\endtitle
\author Meng Chen
\endauthor
\leftheadtext{{\bf M. Chen}}
\rightheadtext{{\bf Pluricanonical maps for threefolds}}
\address \ \ \
\newline Department of Applied Mathematics, Tongji University, Shanghai
200092, China
\newline E-mail: chenmb\@online.sh.cn
\endaddress
\thanks The author was partially supported by the National Natural Science Foundation of China.
\endthanks
\endtopmatter
%%%%%%%%%%%%%%%%%%%%%%%%%%%%
\document

\leftline{\bf 1. Introduction}
\smallskip

This paper is a continuation of \cite{4, 9, 13}. To classify algebraic
varieties is one of the goals
in algebraic geometry. One way to study a given variety is to understand the
behavior
of its pluricanonical maps. The objects concerned here are complex
projective 3-folds
of general type over ${\Bbb C}$.
Let $X$ be such an object and denote by $\phi_m$ the m-th pluricanonical map
of
$X$, which is the rational map associated with the $m$-canonical system
$|mK_X|$.
The very natural question is when $|mK_X|$ gives a birational map, a
generically
finite map, $\cdots$, etc. According to \cite{2, 4, 9, 12, 13}, one has the
following
\proclaim{Theorem 0} Let $X$ be a complex projective 3-fold of general type
with
the canonical index $r$. Then

(i) when $r=1$, $\fei{m}$ is a birational morphism onto its image for $m\ge
6$;

(ii) when $r\ge 2$, $\fei{m}$ is a birational map onto its image for $m\ge
4r+5$.
\endproclaim

In this paper, we give our results on the generic finiteness of $\fei{m}$.
By a delicate use of the Kawamata-Viehweg vanishing theorem, we reduce the
problem
to a parallel one for adjoint systems on some smooth surface. Reider's
results as well
as other theorems on surfaces make it possible for us to go on a detailed
argument.

\proclaim{Theorem 1} Let $X$ be a projective 3-fold of general type with the
canonical
index $r\ge 2$. Then $\phi_m$ is generically finite for $m\ge m(r)$, where
$m(r)$ is a function as follows:

$m(2)=11$;

$m(r)=2r+8$, for $3\le r\le 5$;

$m(r)=2r+6$, for $r\ge 6$.
\endproclaim

\proclaim{Theorem 2} Let $X$ be a projective minimal Gorenstein 3-fold of
general type. Then

(1) $\phi_5$ is birational except for  some 3-folds with $K_X^3=2$ and
$p_g(X)\le 2$;
$\phi_5$ is generically finite of degree $\le 8$.

(2) $\phi_4$ is birational if $K_X^3>2$ and $\dim\phi_1(X)=3$; $\phi_4$ is
generically finite except for some 3-folds with $K_X^3=2$, $p_g(X)\le 1$ and
$\chi(\Co{X})=-1$.

(3) $\phi_3$ is generically finite if $p_g(X)\ge 39$.
\endproclaim

For a nonsingular projective minimal 3-fold $X$ of general type,  Benveniste
(\cite{2}) proved that $\dim\fei{m}(X)\ge 2$ for $m\ge 4$, i.e. $|4K_X|$ can
not be composed of a pencil. Recently, it has been proved (\cite{5}) that
$|3K_X|$ also can not be composed of a pencil.
(Actually, the method is also effective for Gorenstein 3-folds of general
type.)
Thus it is interesting whether $|2K_X|$ can be composed of a pencil
and like what a bicanonical pencil behaves.
So in section 4, we study the bicanonical pencil of a Gorenstein 3-fold of
general type.
According to the 3-dimensional MMP, we can suppose that $X$  is a minimal locally
factorial Gorenstein 3-fold of general type. Take a birational modification
$\pi: X'\lrw X$ such that $X'$ is smooth, $|\pi^*(2K_X)|$ gives a morphism
and
$\pi^*(2K_X)$ has supports with only normal crossings. This is possible
because of Hironaka's big theorem. Let $W:=\overline{\fei{2}(X)}$ and take
the
Stein factorization $$\phi_2\circ\pi: X'\overset f\to\lrw C\overset s\to\lrw
W.$$
Then $f$ is a fibration onto the nonsingular curve $C$, we call $f$ a {\it
derived
fibration} of $\fei{2}$. Denote by $F$ a general fibre of $f$. Then $F$ is a
nonsingular surface of general type by virtue of the Bertini theorem. Also
set
$b:=g(C)$, the geometric genus of $C$. {}From \cite{7}, we know that
$0\le b\le 1$. We shall prove the following
\proclaim{Theorem 3} Let $X$ be a projective minimal Gorenstein 3-fold of
general
type and suppose that $|2K_X|$ is composed of a pencil. Let $f$ be the
derived fibration
of $\fei{2}$ and $F$ be a general fibre of $f$. Then we have
$p_g(F)=1$ and $K_{F_0}^2\le 3$, where $F_0$ is the minimal model of $F$.
\endproclaim

As an application of our method, we shall present a corollary on surfaces of
general type
which somewhat simplifies  Xiao's theorem for the bicanonical finiteness.
\bigskip

\leftline{\bf 2. Proof of Theorem 1}
\smallskip

Throughout our argument, the Kawamata-Viehweg vanishing theorem is always
employed
as a much more effective tool. We use it in the following form.

\proclaim{K-V Vanishing Theorem} (\cite{10} or \cite{17}) Let $X$ be a
nonsingular
complete variety, $D\in\Div(X)\otimes{\Bbb Q}$. Assume the following two
conditions:

(1) $D$ is nef and big;

(2) the fractional part of $D$ has the support with only normal crossings.

\noindent Then $H^i(X,\Co{X}(\roundup{D}+K_X))=0$ for $i>0$, where
$\roundup{D}$
is the round-up of $D$, i.e. the minimum integral divisor with
$\roundup{D}-D\ge 0$.
\endproclaim

\proclaim{Lemma 2.1} (Corollary 2 of \cite{16}) Let $S$ be a nonsingular
algebraic surface,
$L$ be a nef divisor on $S$, $L^2\ge 10$ and let $\phi$ be a map defined by
$|L+K_S|$. If $\phi$ is not birational, then $S$ contains a base point free
pencil $E'$ with $L\cdot E'=1$ or $L\cdot E'=2$.
\endproclaim

\proclaim{Lemma 2.2} Let $X$ be a nonsingular variety of dimension $n$,
$D\in\Div (X)\otimes{\Bbb Q}$ be a ${\Bbb Q}$-divisor on $X$. Then we have
the following:

(i) if $S$ is a smooth irreducible divisor on $X$, then
$\roundup{D}|_S\ge \roundup{D|_S}$;

(ii) if $\pi: X'\lrw X$ is a birational morphism, then
$\pi^*(\roundup{D})\ge \roundup{\pi^*(D)}$.
\endproclaim
\demo{Proof}
We can write $D$ as $G+\sum_{i=1}^t a_iE_i$, where $G$ is a divisor,
the $E_i$ are effective divisors for each $i$ and $0<a_i<1$, $\forall\ i$.
So we only have to prove the lemma for effective ${\Bbb Q}$-divisors.
That is easy to check.
\qed\enddemo

\proclaim{Lemma 2.3} (Lemma 2.3 of \cite{9}) Let $X$ be a minimal threefold
of general type with canonical index $r$. Then we have the plurigenus
formula
$$\align
&h^0(X,\omega_X^{[mr+s]})\\
&=\frac{1}{12}(mr+s)(mr+s-1)(2mr+2s-1)(K_X^3)+am+c_s
\endalign$$
for $0\le s<r$, $mr+s\ge 2$, where $a$ is a constant and $c_s$ is a constant
only relating to $s$.
\endproclaim

\remark{\smc Definition 2.4} Let $X$ be a nonsingular projective variety of
dimension
$\ge 2$. Suppose $|M|$ is a base-point-free system on $X$, {\it a general
irreducible element $S$ of} $|M|$ means the following:

(i) if $\dim\Fi{M}(X)\ge 2$, then $S$ is just a general member of $|M|$;

(ii) if $\dim\Fi{M}(X)=1$, taking the Stein factorization of $\Fi{M}$, then
we
obtain a fibration $f:X\lrw C$ onto a curve $C$. We mean $S$ a general fibre
of $f$.
\endremark

\proclaim{Proposition 2.5} (Lemma 3.2 of \cite{9})
Let $X$ be a minimal threefold of general type with canonical index $r\ge
2$.
Then $\dim\fei{mr+s}(X)\ge 2$ in one of the following cases:

(i) $r=2$ and $m\ge 3$;

(ii) $r=3$ and $m\ge 2$;

(iii) $r=4,5$, $0\le s\le 2$ and $m\ge 2$; $r=4,5$, $s\ge 3$ and $m\ge 1$;

(iv) $r\ge 6$, $0\le s\le 1$ and $m\ge 2$; $r\ge 6$, $s\ge 2$ and
$m\ge 1$.
\endproclaim

Now we modify Proposition 2.5 by virtue of Hanamura's method in order to
prove our Theorem 1. The proof is due to Hamamura (\cite{9}).

\proclaim{Proposition 2.6}
Let $X$ be a minimal threefold of general type with canonical index $r\ge
2$.
Then $h^0(\omega_X^{[mr+s]})\ge 3$ in one of the following cases:

(i) $r=2$ and $m\ge 2$;

(ii) $r\ge 3$, $s=0,1$ and $m\ge 2$; $r\ge 3$, $s\ge 2$ and $m\ge 1$.
\endproclaim
\demo{Proof}
{}From Lemma 2.3, we can put
$$P(mr+s)=\frac{1}{12}(mr+s)(mr+s-1)(2mr+2s-1)(K_X^3)+am+c_s \tag 2.1$$
where $a$ and $c_s$ are constants for $0\le s<r$. We consider the right handside
of (2.1) as a polynomial in $m$ and denote it by $P_s(m)$. Let $Q_s(m)$ be
the first term of $P_s(m)$. We have
$$P_s(m)=Q_s(m)+am+c_s.$$
We see that, for $ m\ge 1$ or $m=0$ and $s\ge 2$,
$$P_s(m)\ge 0. \tag 2.2$$
By Koll\'ar's result (\cite{11}) that the $\omega_X^{[mr+s]}$ are Cohen-Macaulay,
using the Grothendieck duality, one can see that, for $m\le -1$,
$$P_s(m)\le 0.\tag 2.3$$
Now we want to estimate both $a$ and $c_s$. For any $r$ and $s$, by (2.2)
and (2.3), we have
$$Q_s(1)+a+c_s\ge 0\tag 2.4$$
$$-Q_s(-1)+a-c_s\ge 0.\tag 2.5$$

Which induces
$$\align
a&\ge\frac{1}{2}\biggl\{Q_s(-1)-Q_s(1)\biggr\}\tag 2.6\\
&=-\frac{1}{12}\biggl\{2r^2+(6s^2-6s+1)\biggr\}(rK_X^3).
\endalign$$

When $r\ge 3$ and $s\ge 2$, we have
$$Q_s(0)+c_s\ge 0.\tag 2.7$$
By (2.5) and (2.7), we get
$$\align
a&\ge -Q_s(0)+Q_s(-1)\tag 2.8\\
&=\frac{1}{12}\biggl\{-2r^2+(6s-3)r-(6s^2-6s+1)\biggr\}(rK_X^3).
\endalign$$

Explicitly, we have
$$
a\ge\frac{1}{12}\biggl\{-\frac{1}{2}r^2+\frac{1}{2}\biggr\}(rK_X^3)
\ \text{if}\ r\ \text{is odd}
\tag 2.9$$
$$
a\ge\frac{1}{12}\biggl\{-\frac{1}{2}r^2-1\biggr\}(rK_X^3)
\ \text{if}\ r\ \text{is even}.
\tag 2.10$$
Now we can calculate the $P(mr+s)$ case by case.

Case 1. $r\ge 3$ and $s\ge 2$.

When $r$ is odd, from (2.7) and (2.9), we have
$$\align
P(mr+s)&\ge
Q_s(m)-\frac{1}{12}m\biggl(\frac{1}{2}r^2-\frac{1}{2}\biggr)(rK_X^3)-Q_s(0)\\
&=\frac{1}{12}\biggl\{(mr+s)(mr+s-1)(2mr+2s-1)+m\biggl(-\frac{1}{2}r^3+
\frac{1}{2}r\biggr)\\
&\ \ \ \ -s(s-1)(2s-1)\biggr\}(K_X^3).
\endalign$$
We get $P(mr+s)\ge 7$ for $m\ge 1$.

When $r$ is even, from (2.7) and (2.10), we have
$$\align
P(mr+s)&\ge
Q_s(m)-\frac{1}{12}m\biggl(\frac{1}{2}r^2+1\biggr)(rK_X^3)-Q_s(0)\\
&=\frac{1}{12}\biggl\{2r^2m^3+(6s-3)rm^2+(6s^2-6s-\frac{1}{2}r^2)m\biggr\}(r
K_X^3).
\endalign$$
We get $P(mr+s)\ge 5$ for $m\ge 1$.

Case 2. $s=1$.

{}From (2.4) and (2.5), we have
$$P(mr+1)\ge\frac{1}{12}r(m^2-1)(2rm+3)(rK_X^3).$$
We get $P(mr+1)\ge 6$ for $m\ge 2$.

Case 3. $s=0$.

By (2.4) and (2.5), we have
$$P(mr)\ge\frac{1}{12}r(m^2-1)(2rm-3)(rK_X^3).$$
We get $P(mr)\ge 3$ for $m\ge 2$. Thus we complete the proof. \qed
\enddemo

In what follows we can get an improved version of Hanamura's theorem.

\proclaim{Theorem 2.7} Let $X$ be a projective threefold of general type
with the canonical index $r\ge 2$. Then $\fei{m}$ is birational onto its
image
for $m\ge 4r+3$.
\endproclaim
\demo{Proof} We can suppose that $X$ is a minimal 3-fold.
For any  $m_1\ge r+2$,  take some blowing-ups
$\pi: X'\lrw X$ according to Hironaka such that $X'$ is nonsingular and
that the movable part of $|m_1K_{X'}|$ defines a morphism. Denote by $|M|$ the moving part of $|m_1K_{X'}|$
and by $S$ a general irreducible element of $|M|$. Then $S$ is a nonsingular
projective
surface of general type by the Bertini theorem. On $X'$, we consider the
system
$|K_{X'}+3\pi^*(rK_X)+S|$. Because $K_{X'}+3\pi^*(rK_X)$ is effective by Proposition 2.6, so
the system can distinguish general irreducible elements of $|M|$. On the
other hand,
the vanishing theorem gives
$$|K_{X'}+3\pi^*(rK_X)+S||_S=|K_S+3L|,$$
where $L:=\pi^*(rK_X)|_S$ is a nef and big divisor on $S$ and $L^2\ge 2$.
Reider's
result tells that the right system gives a birational map, so does
$|K_{X'}+3\pi^*(rK_X)+S|.$
Thus $\phi_m$ is birational for $m\ge 4r+3$.
\qed\enddemo

\demo{Proof of Theorem 1}
We can suppose that $X$ is a minimal model.
If $r=2$, then $\phi_m$ is birational for $m\ge 11$ according to Theorem
2.7. From now on, we assume $r\ge 3$ and define
$$m_2=\cases r+3, &\text{for $3\le r\le 5$}\\
  r+2, &\text{for $r\ge 6$.}\endcases$$
Take some blowing-ups $\pi:X'\lrw X$ such that $X'$ is nonsingular,
$|m_2K_{X'}|$ defines a morphism and  the fractional part of $\pi^*(K_X)$
has
supports with only normal crossings. Denote by $|M_2|$ the moving part of
$|m_2K_{X'}|$ and by $S_2$ a general irreducible element of $|M_2|$.
For any $t\in {\Bbb Z}_{>0}$, we consider the system
$$|K_{X'}+\roundup{(t+m_2)\pi^*(K_X)}+S_2|,$$
which is a sub-system of $|(t+2m_2+1)K_{X'}|$.
Because $K_{X'}+\roundup{(t+m_2)\pi^*(K_X)}$ is effective by Proposition 2.6,
so the system can
distinguish general irreducible elements of $|M_2|$. On the other hand,
the K-V vanishing theorem tells that
$$\align
&|K_{X'}+\roundup{(t+m_2)\pi^*(K_X)}+S_2||_{S_2}\\
&=|G+L|,
\endalign$$
where $G:=\{K_{X'}+\roundup{(t+m_2)\pi^*(K_X)}\}|_{S_2}$ is effective and
$L:=S_2|_{S_2}$. We can see that
$$G+L\ge K_{S_2}+\roundup{t\pi^*(K_X)}|_{S_2}+L.$$
{}From Proposition 2.5, we have $h^0(S_2, L)\ge 2$.
Modulo blowing-ups, actually we can suppose that $|L|$ is free from base points.
Let $C$ be a general irreducible element of $|L|$. It is obvious that
$ |G+L|$ can distinguish gereral irreducible elements of $|L|$. On the other
hand,
the K-V vanishing theorem gives
$$|K_{S_2}+\roundup{t\pi^*(K_X)|_{S_2}}+C||_C=|K_C+D|,$$
where $D:=\roundup{t\pi^*(K_X)|_{S_2}}|_C$ is a divisor of positive degree.
Because
$C$ is a curve of genus $\ge 2$, so $h^0(C, K_C+D)\ge 2$ and $|K_C+D|$ gives
a
finite map. Thus we have $\dim\Fi{G+L}(C)=1$. Therefore $\phi_m$
is generically finite for $m\ge 2m_2+2$, which completes the proof.
\qed\enddemo
\bigskip

\leftline{\bf 3. On Gorenstein 3-folds of general type}
\smallskip

For a minimal threefold $X$ of general type with canonical index $1$, we can
find certain birational modifications $f:X'\lrw X$ according to \cite{15}
such that $c_2(X')\cdot\Delta=0$, where $\Delta$ is the ramification
divisor of $f$.
Then we can get the same plurigenus formula as that for a nonsingular
minimal
threefold, i.e.
$$p(n):=h^0(X,
\Co{X}(nK_X))=(2n-1)\bigl\{n(n-1)K_X^3/12-\chi(\Co{X})\bigr\},$$
for $n\ge 2$.
On the other hand, the Miyaoka-Yau inequality (\cite{14}) shows that
$\chi(\Co{X})<0$.
{}From \cite{4} or \cite{12}, we know that $\fei{m}$ is  birational for
$m\ge 6$.

\proclaim{Theorem 3.1} Let $X$ be a projective minimal Gorenstein 3-fold of
general type. Then

(1) $\fei{5}$ is  birational if either $K_X^3>2$ (Ein-Lazarsfeld-Lee) or
$p_g(X)>2$.

(2) When $p_g(X)=2$, then $\fei{5}$ is birational except for some 3-folds
with
$q(X)=h^2(\Co{X})=0$, and $|K_X|$ composed with a
rational pencil of surfaces of general type with $(K^2,p_g)=(1,2)$. In this
situation, $\fei{5}$ is generically finite of degree $2$.

(3) $\fei{5}$ is birational if $\dim\phi_2(X)=1$.
\endproclaim
\demo{Proof} This is the main theorem in \cite{7}. Though the objects
considered
there are nonsingular minimal 3-folds, the method is also effective for all
Gorenstein 3-folds of general type.
\qed\enddemo

\remark{\smc Definition 3.2} Let $X$ be a projective minimal Gorenstein
3-fold of general type. Suppose $\dim\fei{i}(X)\ge 2$ and set
$iK_X\simlin M_i+Z_i$, where $M_i$ is the moving part and $Z_i$ the fixed
one
for any integer $i$.
We define $\delta_i(X):=K_X^2\cdot M_i$.
\endremark

\proclaim{Proposition 3.3} Let $X$ be a projective minimal Gorenstein 3-fold
of general type.
Suppose $|2K_X|$ is not composed of a pencil and $K_X^3>2$. Then
$\delta_2(X)\ge 3$.
\endproclaim
\demo{Proof} We have $\delta_2(X)\ge 2$ by Proposition 2.2 of \cite{4}.
Take a birational modification $f:X'\lrw X$ such that $|2f^*(K_X)|$ defines
a morphism. Set $2f^*(K_X)\simlin M+Z$, where $M$ is the moving part and $Z$
the
fixed one. A general member $S\in|M|$ is an irreducible nonsingular
projective surface of general
type. Denote $L:=f^*(K_X)|_S$. If $L^2={f^*(K_X)}^2\cdot S=\delta_2(X)=2$,
then we have
$$4=2{f^*(K_X)}^2\cdot S=f^*(K_X)\cdot S^2+f^*(K_X)\cdot S\cdot Z.$$
Noting that $S$ is nef and $S\not\approx 0$, we have $f^*(K_X)\cdot S^2\ge 1$.
Therefore four cases occur as follows:

(i) $f^*(K_X)\cdot S^2=4$, $f^*(K_X)\cdot S\cdot Z=0;$

(ii) $f^*(K_X)\cdot S^2=3$, $f^*(K_X)\cdot S\cdot Z=1;$

(iii) $f^*(K_X)\cdot S^2=2$, $f^*(K_X)\cdot S\cdot Z=2;$

(iv) $f^*(K_X)\cdot S^2=1$, $f^*(K_X)\cdot S\cdot Z=3.$

We also have
$$\align
2K_X^3&=2f^*(K_X)^3=f^*(K_X)^2\cdot S+f^*(K_X)^2\cdot Z \tag 3.1\\
      &=2+\frac{1}{2}f^*(K_X)\cdot Z(S+Z)\\
      &=2+\frac{1}{2}f^*(K_X)\cdot S\cdot Z+\frac{1}{2}f^*(K_X)\cdot Z^2.
\endalign$$

\ \ \ Case (i). Noting that $f^*(K_X)$ is nef and big, we see that
$mf^*(K_X)$ is
linearly equivalent to a nonsingular projective surface of general type
according
to Kawamata for sufficiently large integer $m$. Then $S|_{mf^*(K_X)}$ is
nef and big and,
by the Hodge Index Theorem, we have $f^*(K_X)\cdot Z^2\le 0$. Thus (3.1)
is false and this case does not occur.

Case (ii). We have $f^*(K_X)\cdot S(S-3Z)=0$, then $f^*(K_X)(S-3Z)^2\le 0$,
which derives $f^*(K_X)\cdot Z^2\le \frac{1}{3}$, i.e. $f^*(K_X)\cdot Z^2\le
0$.
(3.1) is also false.

Case (iii). $f^*(K_X)\cdot S(S-Z)=0$ induces $f^*(K_X)\cdot Z^2\le 2$, then
(3.1) becomes $K_X^3\le 2$. Thus $K_X^3=2$. Actually, in this case,
$f^*(K_X)\cdot (S-Z)\simnum 0$ (as 1-cycle).

Case (iv). $f^*(K_X)\cdot (3S-Z)^2\le 0$ induces $f^*(K_X)\cdot Z^2\le 9$.
And (3.1) becomes $K_X^3\le 4$. If $K_X^3=4$, we see that $f^*(K_X)\cdot
(3S-Z)
\simnum 0$ as 1-cycle. Now we set $f^*(M_2)=S+E$. Then $Z=f^*(Z_2)+E$.
Obviously, we have $f_*(S)=M_2$ and $f_*(Z)=Z_2$.
{}From $f^*(M_2)\cdot f^*(K_X)\cdot (3S-Z)=0$, we get $3K_X\cdot
M_2^2=K_X\cdot M_2
\cdot Z_2$. Then $4=2K_X^2\cdot M_2=K_X\cdot M_2^2+K_X\cdot M_2\cdot Z_2=
4K_X\cdot M_2^2$, i.e. $K_X\cdot M_2^2=1$. Which derives a contradiction,
because $K_X\cdot M_2^2$ is even. Thus $K_X^3=2$.
\qed\enddemo

\proclaim{Proposition 3.4} Let $X$ be a projective minimal Gorenstein 3-fold
of general type.
Suppose $K_X^3>2$ and $\dim\fei{1}(X)\ge 2$. Then $\delta_1(X)\ge 3$.
\endproclaim
\demo{Proof} As in the proof of the previous proposition, we first take
a modification $f: X'\lrw X$. Set $f^*(K_X)\simlin M+Z$, where $M$ is
the moving part. A general member $S\in|M|$ is a nonsingular projective
surface
of general type. Also denote $L:=f^*(K_X)|_S$. Then $L^2=\delta_1(X)\ge 2$
according to Proposition 2.1 of \cite{7}. If $L^2=2$, then we have
$$2=f^*(K_X)^2\cdot S=f^*(K_X)\cdot S^2+f^*(K_X)\cdot S\cdot Z.$$
We also have
$$\align
K_X^3&=f^*(K_X)^2\cdot S+f^*(K_X)^2\cdot Z \tag 3.2\\
     &=2+f^*(K_X)\cdot S\cdot Z+f^*(K_X)\cdot Z^2.
\endalign$$
Similarly, $f^*(K_X)\cdot S^2\ge 1$. If $f^*(K_X)\cdot S^2=2$ and
$f^*(K_X)\cdot S\cdot Z=0$, then, by the Hodge Index Theorem, $f^*(K_X)\cdot
Z^2\le 0$.
Then (3.2) becomes $K_X^3\le 2$, which says $K_X^3=2$.
If  $f^*(K_X)\cdot S^2=f^*(K_X)\cdot S\cdot Z=1$, $f^*(K_X)\cdot S\cdot
(S-Z)=0$
induces $f^*(K_X)\cdot Z^2\le 1$. By (3.2), we get $K_X^3\le 4$.
If $K_X^3=4$, then we can see $f^*(K_X)\cdot(S-Z)\simnum 0$. 
By the same argument as in the case (iv) of the proof of Proposition 3.3,
we have
$f^*(M_1)\cdot f^*(K_X)\cdot(S-Z)=0$, i.e. $K_X\cdot M_1^2=K_X\cdot M_1\cdot
Z_1$.
We have $2=K_X^2\cdot M_1=K_X\cdot M_1^2+K_X\cdot M_1\cdot Z_1=2K_X\cdot
M_1^2$.
Therefore $K_X\cdot M_1^2=1$, which is impossible. Thus $K_X^3=2$.
\qed
\enddemo

\proclaim{Theorem 3.5} Let $X$ be a projective minimal Gorenstein 3-fold of
general type.
Then $\fei{5}$ is generically finite
of degree $\le 8$. If $\deg(\fei{5})>2$, then $K_X^3=2$, $\chi(\Co{X})=-1$
and $p_g(X)=0,1$.
\endproclaim
\demo{Proof}
According to Theorem 3.1, we only have to study the case when $|2K_X|$ is
not
composed of a pencil. Take a modification $f:X'\lrw X$ according to Hironaka
such that $|2f^*(K_X)|$ defines a morphism.
Set $2f^*(K_X)\simlin M+Z$, where $M$ is the moving part and
$Z$ the fixed one. A general member $S\in|M|$ is a nonsingular projective
surface
of general type by the Bertini Theorem. We have
$$|K_{X'}+2f^*(K_X)+S|\subset|5K_{X'}|.$$
Because $K_{X'}+2f^*(K_X)$ is effective, the left system can distinguish
general
members of $|M|$.
Denote $L:=f^*(K_X)|_S$, using the long exact sequence and the
vanishing theorem, we have
$$|K_{X'}+2f^*(K_X)+S|\vert_S=|K_S+2L|.$$
Obviously, $K_S+2L=G+H$, where $G:=(K_{X'}+2f^*(K_X))|_S$ is effective and
$H:=S|_S$.
Note that $h^0(S,\Co{S}(2L))\ge h^0(S,H)\ge P(2)-1\ge 3$. We have two cases.

Case 1. $|H|$ is composed of a pencil. Taking a birational modification to $S$
if
necessary, we can suppose $|H|$ is free from base points.
Denote $H\simlin
\sum_{i=1}^a C_i+E,$ where  $E$ is the fixed part. In general position,
$\sum_{i=1}^a C_i$ can be a disjoint union of nonsingular curves in a
family.
We have $a\ge 2$. Thus $L\simnum \frac{a}{2}C+E_0$, where
$E_0\ge\frac{1}{2}E$
is an effective ${\Bbb Q}$-divisor.
If $p_g(S)=0$, then $q(S)=0$ and then we can see by the long exact sequence
that $|K_S+H|$ can distinguish $C_i$'s and that
$|K_S+\sum_{i=1}^aC_i|\vert_{C_i}=|K_{C_i}|$,
which means $|K_S+2L|$ gives at worst a generically finite map of degree $2$
and
so does $\phi_5$. If
$p_g(S)>0$, it is obvious that $|K_S+2L|$ can distinguish $C_i$'s. For a
general
curve $C$ which is algebraically equivalent to $C_i$, we consider the
${\Bbb Q}$-divisor $G:=K_S+2L-\frac{1}{2}\sum_{i=3}^aC_i-E_0$.
We have $\roundup{G}\le K_S+2L$. On the other hand, $G-C-K_S$ is nef and big,
thus by the K-V vanishing we have
$|\roundup{G}|\vert_C=|K_C+\roundup{E_0}|_C|$. Because $\roundup{E_0}|_C$ is effective, 
$\Fi{K_S+2L}$ is at worst a generically finite map of degree 2 and so is
$\fei{5}$ of $X$.

Case 2. $|H|$ is not composed of a pencil, so neither is $|2L|$. Similarly,
we can suppose $|2L|$
is base point free. If $p_g(S)=0$, we can use a parallel discussion to that
of Case 1 to see that $\fei{5}$ is at worst a generically finite map of
degree
$2$. If $p_g(S)>0$, then $\Fi{K_S+2L}$ is obviously generically finite.
We know that $L^2\ge 2$ from Proposition 2.2 of \cite{4}. If
$\Fi{K_S+2L}$ is not birational and $L^2\ge 3$,
then according to Lemma 2.1, there is a free pencil on $S$
with a general member $C$ such that $C^2=0$ and $L\cdot C=1$. Since
$\dim\Fi{2L}(C)=1$, then $h^0(2L|_C)\ge 2$ and then, by the Clifford
theorem,
we see that $C$ is a curve of genus 2 and $2L|_C\simlin K_C$. Finally we can
see that
$|2L|\vert_C=|K_C|$. Therefore $\Fi{K_S+2L}$ is a generically finite map of
degree  2. Therefore $\fei{5}$ is generically finite with $\deg(\fei{5})\le
2$.
If $L^2=2$, then $K_X^3=2$ by the proof of Proposition 3.3. On the surface $S$, set
$2L\simlin C_1+E_1$, where $C_1$ is the moving part. We easily get
$$8=(2L)^2\ge C_1^2\ge d(h^0(2L)-2)\ge d(P(2)-3).$$
Therefore we have
$$d\le\frac{8}{P(2)-3}=\frac{8}{-3\chi(\Co{X})-2}.$$
If $d>2$, then $\chi(\Co{X})=-1$. \qed
\enddemo

For the 4-canonical map of $X$, it is obvious that $\fei{4}$ is not
birational if
$X$ admits a  pencil of  surfaces of general type with $(K^2,p_g)=(1,2)$.
Therefore it is pessimistic for us to obtain an effective sufficient
condition
for the birationality of $\fei{4}$. We have a partial result as follows.

\proclaim{Theorem 3.6} Let $X$ be a projective minimal Gorenstein 3-fold of
general type.
Suppose $K_X^3>2$ and $\dim\fei{1}(X)=3$. Then $\fei{4}$ is
a birational map onto its image.
\endproclaim
\demo{Proof}
Take a birational modification $f:X'\lrw X$ such that the movable part of $|f^*(K_X)|$ is base point free.
Set $f^*(K_X)\simlin S+Z$, where $S$ is the moving part and $Z$ the fixed
one.
A general member $S$ is a nonsingular projective surface of general type. We
have
$|K_{X'}+2f^*(K_X)+S|\subset|4K_{X'}|$. Using the vanishing theorem, we have
$$|K_{X'}+2f^*(K_X)+S|\vert_S=|K_S+2L|,$$
where $L:=f^*(K_X)|_S$ is a nef and big divisor on $S$. By Proposition 3.4,
we see that $L^2\ge 3$ under the condition $K_X^3>2$. If $\Fi{K_S+2L}$ is
not birational,
then, by Lemma 2.1, there is a free pencil with a general member
$C$ such that $C^2=0$ and $L\cdot C=1$. Because $\dim\Fi{L}(S)=2$, $h^0(C,
\Co{C}(L|_C))\ge 2$. Therefore, by the Clifford theorem, we see that
$\deg(L|_C)\ge 2h^0(L|_C)-2\ge 2$. This is a contradiction. Therefore
$\Fi{K_S+2L}$ is birational and so is $\fei{4}$.\qed
\enddemo

\remark{\smc Example 3.7} We give an example which shows that $\fei{4}$ is
not birational
when $K_X^3=2$ and $\dim\fei{1}(X)=3$. On ${\Bbb P}^3({\Bbb C})$, take a
smooth
hypersurface $S$ of degree 10, $S\simlin 10H$. Let $X$ be a double cover of
${\Bbb P}^3$ with branch locus along $S$. Then
$X$ is a nonsingular canonical model,
$K_X^3=2$ and $p_g(X)=4$ and
$\fei{1}$ is a finite morphism onto ${\Bbb P}^3$ of degree $2$.
One can easily check that $\fei{4}$ is also a finite morphism of degree $2$.
\endremark

\proclaim{Theorem 3.8} Let $X$ be a projective minimal Gorenstein 3-fold of
general type.
Then $\fei{4}$ is generically finite when $p_g(X)\ge 2$ or when $K_X^3>2$
or when $\chi(\Co{X})\not= -1$.
\endproclaim
\demo{Proof}
Part I: $p_g(X)\ge 2$.

First we make a modification $f:X'\lrw X$ such that the movable part of 
$|f^*(K_X)|$ is free from base points and that $f^*(K_X)$ has support with
only normal
crossings. Set $f^*(K_X)\simlin M+Z$, where $M$ is the moving part and $Z$
the fixed  one.

If $\dim\fei{1}(X)=2$, then a general member $S\in|M|$ is a nonsingular
projective
surface of general type. We have
$$|K_{X'}+2f^*(K_X)+S|\subset|4K_{X'}|.$$
Using the vanishing theorem, we have $|K_{X'}+2f^*(K_X)+S||_S=|K_S+2L|$,
where
$L:=f^*(K_X)|_S$ is nef and big effective divisor on $S$. We have
$h^0(S, L)\ge 2$. Noting that $p_g(S)>0$ in this case. And if $|L|$ is not
composed
of a pencil, then neither is $|K_S+2L|$. If $|L|$ is composed of a pencil,
taking a
modification if possible, we can suppose that the movable part of $|L|$ is free from base points.
Set $L\simlin \sum C_i+Z_0$, we can see $|K_S+L+\sum C_i||_{C_i}=
|K_{C_i}+D|$, where $D:=L|_{C_i}$ is effective. We easily see that
$\Fi{K_S+2L}$ is at worst generically finite of degree $\le 2$ and so is
$\fei{4}$.

If $\dim\fei{1}(X)=1$, then $M\simnum aF$, where $F$ is a nonsingular
projective
surface of general type. $M_1\simnum aF_0$, where $F_0=f_*(F)$ is
irreducible
on $X$. If $K_X\cdot F_0^2=0$, then, by Lemma 2.3 of \cite{7}, we have
$\Co{F}(f^*(K_X)|_F)\cong\Co{F}(\pi^*(K_0))$, where $\pi$ is the contraction
map
onto the minimal model and $K_0$ is the canonical divisor of the minimal
model of
$F$. Obviously, $|K_{X'}+2f^*(K_X)+M|$ can distinguish general members of
$|M|$.
Moreover $|K_{X'}+2f^*(K_X)+M||_F=|K_F+2\pi^*(K_0)|$, the right system gives
a
generically finite map and so does $\phi_4$.
If $K_X\cdot F_0^2>0$,  then
$$L^2=f^*(K_X)^2\cdot F=K_X^2\cdot F_0\ge K_X\cdot F_0^2\ge 2.$$
It is sufficient to show that $|K_F+2L|$ gives a generically finite map.
We have $K_F+2L\ge 3L$. If $|3L|$ is not composed of a pencil, then neither
is
$|K_F+2L|$. If $|3L|$ is composed of a pencil, we claim that $h^0(F, 3L)\ge
3$.
In fact, we have $|K_{X'}+f^*(K_X)+F||_F=|K_F+L|$ and  $h^0(F, K_F+L)\ge 3$.
Considering the natural map $H^0(X', 3K_{X'})\overset\alpha\to\lrw
H^0(F,3K_F)$,
because $K_{X'}+f^*(K_X)+F\le 3K_{X'}$, we see that
$\dim_{\Bbb C}(\text{Im}(\alpha))\ge h^0(K_F+L)\ge 3$. Similarly,
considering
another natural map $H^0(X', 3f^*(K_X))\overset\beta\to\lrw H^0(F, 3L)$,
we have
$$h^0(3L)\ge \dim_{\Bbb C}(\text{Im}(\beta))=\dim_{\Bbb
C}(\text{Im}(\alpha))
\ge 3.$$
Now we can write $3L\simlin\sum_{i=1}^t\overline{C_i}+E_0$, where $E_0$ is the fixed part, $t\ge 2$ and the $\overline{C_i}$ are irreducible curves. Denote by $C$ a generic $\overline{C_i}$. Then $2L\simnum \frac{2}{3}tC+\frac{2}{3}E_0$ and thus
$2L-C-\frac{1}{t}E_0$ is a nef and big ${\Bbb Q}$-divisor. Setting
$G:=2L-\frac{1}{t}E_0$, then we have $K_S+\roundup{G}\le K_S+2L$. On the other hand,
the K-V vanishing gives $|K_S+\roundup{G}||_C=|K_C+D|,$ where $D$ is a
divisor
of positive degree. Noting that $C$ is a curve of genus $\ge 2$, so we see that
$|K_C+D|$
gives a generically finite map. This means $|K_S+2L|$ gives a generically
finite map.

Part II: $K_X^3>2$ or $\chi(\Co{X})\not=-1$.

We study $\fei{4}$ according to the behavior
of $\fei{2}$. Of course, first we make a modification $f:X'\lrw X$ such that
the movable part of $|2f^*(K_X)|$ is free from base points and that $2f^*(K_X)$ has supports with
only normal crossings. Set $2f^*(K_X)\simlin \overline{M_2}+\overline{Z_2}$,
where $\overline{M_2}$ is the moving part
and $\overline{Z_2}$ the fixed one.

If $\dim\fei{2}(X)=1$, then $\overline{M_2}\simnum a_2F$, where $F$ is a
nonsingular projective
surface of general type. We have
$\Co{F}(f^*(K_X)|_F)\cong\Co{F}(\pi^*(K_0))$
by Lemma 4.2 below in this paper.
Because $K_{X'}+f^*(K_X)$ is effective, $|K_{X'}+f^*(K_X)+\overline{M_2}|$
can
distinguish general $F$. On the other hand,
we have $|K_{X'}+f^*(K_X)+\overline{M_2}||_F=|K_F+\pi^*(K_0)|$.
{}From Theorem 3.1 of \cite{7}, we know that $F$ is not a surface with
$p_g=q=0$. Thus $|K_F+\pi^*(K_0)|$ defines a generically finite map
according to
\cite{19} and so does $\fei{4}$.

If $\dim\fei{2}(X)\ge 2$, then a general member $S\in|\overline{M_2}|$ is a
nonsingular
projective surface of general type. We have
$|K_{X'}+f^*(K_X)+S||_S=|K_S+L|$,
where $L:=f^*(K_X)|_S$. Noting that $K_S\ge L$, then we have $K_S+L\ge 2L$. Under our
assumption, we have $P(2)\ge 5$. Thus $h^0(2L)\ge 4$. We may suppose that the movable part of $|2L|$ is free from base points.
If $|2L|$ is not composed of a pencil, then neither is $|K_S+L|$.
Otherwise we can set $2L\simlin \sum_{i=1}^{b}C_i+E_1$, where
$b\ge 3$ and $E_1$ is the fixed part. We denote by $C$ the general $C_i$.
Because $L-C-\frac{1}{b}E_1$ is nef and big, therefore
$$|K_S+\roundup{L-\frac{1}{b}E_1}|\vert_C=|K_C+D|,$$
where $D$ is a divisor of positive degree.
The right system obviously defines a generically finite map. Thus
$|K_S+L|$ gives a generically finite map and so does
$\fei{4}$.
\qed\enddemo

\proclaim{Theorem 3.9} Let $X$ be a projective minimal Gorenstein 3-fold of
general type.
Then $\fei{3}$ is generically finite when $p_g(X)\ge 39$.
\endproclaim
\demo{Proof}
First we make a modification $f:X'\lrw X$ such that the movable part of 
$|f^*(K_X)|$ is free from base points and that $f^*(K_X)$ has support with
only normal
crossings. Set $f^*(K_X)\simlin M+Z$, where $M$ is the moving part and $Z$
the
fixed  one.

If $\dim\fei{1}(X)\ge 2$, then a general member $S\in|M|$ is a nonsingular
projective
surface of general type. We have $|K_{X'}+f^*(K_X)+S||_S=|K_S+L|$, where
$L:=f^*(K_X)|_S$. When $p_g(X)\ge 4$, $h^0(S,L)\ge 3$. Noting that
$p_g(S)>0$,
if $|L|$ is not composed of a pencil, then nor is $|K_S+L|$. So we may
suppose that
$|L|$ is composed of a pencil and the movable part of this system is free from base points. Set $L\simlin
\sum_{i=1}^{a}C_i+E_0$, where we have $a\ge 2$.
$|K_S+L|$ can distinguish the $C_i$ generically. On the other hand,
$L-C-\frac{1}{a}E_0$ is nef and big, we obtain by the Kawamata-Viehweg
vanishing that
$$|K_S+\roundup{L-\frac{1}{a}E_0}||_C=|K_C+\roundup{\frac{a-1}{a}L}|_C|.$$
The right system defines a generically finite map and so does $\fei{3}$.

If $\dim\fei{1}(X)=1$,
then $M\simnum aF$, where $F$ is a nonsingular projective
surface of general type. Set $F_0=f_*(F)$.
If $K_X\cdot F_0^2=0$, then, by Lemma 2.3 of \cite{7}, we have
$\Co{F}(f^*(K_X)|_F)\cong\Co{F}(\pi^*(K_0))$, where $\pi$ is the contraction
onto
the minimal model and $K_0$ is the canonical divisor of the minimal model of
$F$.
We see that  $|K_{X'}+f^*(K_X)+M||_F=|K_F+\pi^*(K_0)|$.
Because $p_g(F)>0$, the right system defines
a generically finite map and so does $\fei{3}$.
If $K_X\cdot F_0^2>0$, in order to prove the theorem, we have to show
the generic finiteness of $\Fi{K_F+L}$, where $L:=f^*(K_X)|_F$ is effective.
By Theorem 2 of
\cite{6}, we see that $q(F)\ge 3$ when $p_g(X)\ge 39$.
Then $\Fi{K_F}$ is generically finite according to \cite{18}.
Therefore under the
assumption of the theorem, we can obtain the generic finiteness of
$\fei{3}$.
\qed
\enddemo
\bigskip

\leftline{\bf 4. On bicanonical systems}
\smallskip

We suppose that $X$ is a locally factorial Gorenstein minimal 3-fold of general
type and that
$|2K_X|$ be composed of a pencil. Keep the same notations as in section 1 and
let
$\pi:X'\lrw X$ be the birational modification and $f:X'\lrw C$ be the
derived fibration.

\proclaim{Lemma 4.1} Let $X$ be a projective minimal Gorenstein 3-fold of
general type and suppose that $|2K_X|$ is composed of a pencil. Then $q(X)\le 2$ and $p_g(X)\ge 1$.
\endproclaim
\demo{Proof} This is just a generalized version of Corollary 3.1 of
\cite{7}.
Though the objects considered there are nonsingular minimal 3-folds, the
method
is also effective for minimal Gorenstein 3-folds.
\qed\enddemo

\proclaim{Lemma 4.2} Let $X$ be a projective minimal Gorenstein 3-fold of
general
type, $|2K_X|$ be composed of a pencil, $f: X'\lrw C$ be the derived
fibration
of $\phi_2$ and $F$ be a general fibre of $f$. Then
$$\Co{F}(\pi^*(K_X)|_F)\cong\Co{F}(\pi_0^*(K_{F_0})),$$
where $\pi_0: F\lrw F_0$ is the birational contraction onto the minimal
model.
\endproclaim
\demo{Proof} This is just a generalized version of Corollary 9.1 of
\cite{13}.
Though the objects considered there are nonsingular minimal 3-folds, the
method
is also effective for minimal Gorenstein 3-folds.
\qed\enddemo

\proclaim{Lemma 4.3} Under the same assumption as in Lemma 4.2, we have
$K_{F_0}^2\le 3$ and $1\le p_g(F)\le 3$.
\endproclaim
\demo{Proof}
Let $\pi^*(2K_X)\simlin g^*(H_2)+Z_2'$, where $g:=\phi_2\circ\pi$, $Z_2'$ is
the fixed part and $H_2$ is a general hyperplane section of the closure $W$ of the image of $X$ in
${\Bbb P}^{p(2)-1}$. Obviously we have $g^*(H_2)\simnum a_2F$, where $a_2\ge
p(2)-1$.
{}From Lemma 4.2, we have
$$K_{F_0}^2=(\pi^*(K_X)|_F)^2=\pi^*(K_X)^2\cdot F.$$
Let $2K_X\simlin M_2+Z_2$, where $M_2$ is the moving part and $Z_2$ is the
fixed
part. We also have $M_2=\pi_*(g^*(H_2))$. Denote $\overline{F}:=\pi_*(F)$,
then
$M_2\simnum a_2\overline{F}.$ By the projection formula, we get
$$K_X^2\cdot\overline{F}=\pi^*(K_X)^2\cdot F=K_{F_0}^2.$$
Because $K_X$ is nef and big, we have $2K_X^3\ge a_2K_X^2\cdot\overline{F}$.
Thus
$$K_X^2\cdot\overline{F}\le\frac{2}{a_2}K_X^3\le
\frac{4K_X^3}{K_X^3-6\chi(\Co{X})-2}\le\frac{4K_X^3}{K_X^3+4}<4,$$
which means $K_{F_0}^2\le 3$. By Lemma 4.1, the fact that $p_g(X)\ge 1$
induces
$p_g(F)>0$. By the Noether inequality $2p_g(F_0)-4\le K_{F_0}^2$, we see
that
$p_g(F)\le 3$.
\qed\enddemo

\demo{Proof of Theorem 3}
In order to prove Theorem 3, we shall derive a contradiction under the
assumption
that $p_g(F)\ge 2$. Obviously, $|2K_{X'}|$ can distinguish general fibres of
the morphism $\fei{2}\circ \pi$. We consider the system
$|K_{X'}+\pi^*(K_X)|$.
Write $2\pi^*(K_X)\simlin M_2'+Z_2'$, where $M_2'$ is the moving part and
$Z_2'$
is the fixed one.
Set
$Z_2'=Z_v+Z_h$, where $Z_v$ is the vertical part and $Z_h$ is the horizontal
part
with respect to the fibration $f:X'\lrw C$. Noting that $\pi^*(K_X)$ is
effective
by Lemma 4.1, $Z_h$ should be 2-divisible, i.e. $Z_h=2Z_0$, where $Z_0$ is
an
effective divisor. Thus we see that $Z_0$ is just the horizontal part of
$\pi^*(K_X)$. We know that $a_2\ge p(2)-1\ge 3$ and
$$\pi^*(K_X)\simnum \frac{a_2}{2}F+\frac{1}{2}Z_2'.$$
Therefore $\pi^*(K_X)-F-\frac{1}{a_2}Z_2'$ is a nef and big ${\Bbb
Q}$-divisor.
Setting  $G:=\pi^*(K_X)-\frac{1}{a_2}Z_2'$, then we have
$K_{X'}+\roundup{G}\le K_{X'}+\pi^*(K_X)$. By the Kawamata-Viehweg vanishing
theorem,
we see that, for a general fibre $F$,
$$|K_{X'}+\roundup{G}||_F=|K_F+\roundup{G}|_F|\supset|K_F+\roundup{G|_F}|
=|K_F+\roundup{\frac{a_2-2}{a_2}Z_0|_F}|,$$
where $\roundup{\frac{a_2-2}{a_2}Z_0|_F}$ is effective on the surface $F$.
This
means that $\dim\phi_2(F)\ge 1$ under the assumption $p_g(F)\ge 2$ and then
$\dim\phi_2(X)\ge 2$, a contradiction.
\qed\enddemo

The rest of this section is devoted to present an application of our method
to bicanonical
maps of surfaces of general type.

\proclaim{Theorem 4.4} Let $S$ be a minimal algebraic surface of general
type with $p(2)\ge 4$. Then the bicanonical map of $S$ is generically
finite.
\endproclaim
\demo{Proof}
Suppose that $|2K_S|$ is composed of a pencil, we want to derive a
contradiction. Taking a birational modification $\pi: S'\lrw S$ such that
$|2\pi^*(K_S)|$ defines a morphism and denoting $W:=\overline{\phi_2(S)}$, we
obtain the following through the Stein factorization:
$$\phi_2\circ\pi: S'\overset f\to\lrw B\lrw W,$$
where $B$ is a nonsingular curve. Denote by $C$ a general fibre of the
derived
fibration $f$. We can write
$$\pi^*(2K_S)\simlin \sum_{i=1}^{a}C_i+Z,$$
where $a\ge p(2)-1\ge 3$ and $Z$ is the fixed part. Considering the system
$|K_{S'}+\pi^*(K_S)|$, we can see that the system can distinguish general fibres of $\phi_2$. Setting
$G:=\pi^*(K_S)-\frac{1}{a}Z$, we have $K_S+\roundup{G}\le K_S+\pi^*(K_S)$ and
$G-C\simnum \frac{a-2}{a}\pi^*(K_S)$ is nef and big. Thus, by the K-V
vanishing
theorem, we have
$$|K_S+\roundup{G}||_C=|K_C+D|,$$
where $D:=\roundup{G}|_C$ is a divisor of positive degree on the curve $C$.
Because $g(C)\ge 2$, then $h^0(C, K_C+D)\ge 2$. This means that
$|K_S+\pi^*(K_S)|$ gives a generically finite map, a contradiction.
\qed\enddemo

\proclaim{Corollary 4.5} Let $S$ be a minimal algebraic surface of general
type
with $p_g\ge 2$. Then the bicanonical map of $S$ is generically finite.
\endproclaim
\demo{Proof} If $q=0$, then $\chi(\Co{S})\ge 3$ and $p(2)\ge 4$.
If $q>0$, then $K_S^2\ge 2p_g\ge 4$ by \cite{8} and then $p(2)\ge 5$.
The proof is completed by Theorem 4.4.
\qed\enddemo

\proclaim{Corollary 4.6} Let $S$ be a minimal algebraic surface of general
type
with $p(2)=3$. Then $|2K_S|$ is not composed of an irrational pencil.
\endproclaim
\demo{Proof}  This is obvious from the proof of Theorem 4.4. The critical
point
is that we also have $a\ge 3$ in this case.
\qed\enddemo

The remain cases are like the following:

(I) \ \  $K^2=1$, $p_g=1$ and $q=0$;

(II) $K^2=2$ and $p_g=q=0$;

(III) $K^2=2$ and  $p_g=q=1$.

\proclaim{Proposition 4.7} Let $S$ be a minimal algebraic surface of type
(I). Then
the bicanonical map is generically finite.
\endproclaim
\demo{Proof}
Suppose that $|2K_S|$ is composed of a rational pencil. We write
$$2K_S\simlin C_1+C_2+Z,$$
where $Z$ is the fixed part. Denote by $C$ a general member which is
algebrally equivalent to $C_i$. We have $1=K_S^2\ge K_S\cdot C$. On the
other hand,
$K_S\cdot C+C^2\ge 2$, which gives $C^2\ge 1$. Thus $K_S\cdot C=C^2=1$, i.e.
$C$ is a nonsingular curve of genus two. By the index theorem, we see that
$K_S\simnum C$. But from \cite{3}, $\text{Pic}(S)$ is torsion free, then
$K_S\simlin C$. This is impossible because $h^0(S, C)=2$.
\qed\enddemo

\proclaim{Lemma 4.8} (Lemma 8 of \cite{19}) Let $S$ be a surface with finite
$\pi_1$. Then
$$H^1(S, \Co{S}({\Cal E}))=0$$
for any invertible torsion sheaf ${\Cal E}$ on $S$.
\endproclaim

\proclaim{Lemma 4.9} Let $S$ be a minimal surface of type (II) or (III).
Suppose that $|2K_S|$
is composed of a rational pencil. Then the moving part of $|2K_S|$ is a free
pencil of genus two.
\endproclaim
\demo{Proof} We can write $2K_S\simlin C_1+C_2+Z$, where $Z$ is the fixed
part.
Denote by $C$ the general member which is algebraically equivalent to $C_i$.
If $C^2>0$, then $K_S^2\ge K_S\cdot C\ge C^2$. On the other hand, the index
theorem
gives $K_S^2\times C^2\le (K_S\cdot C)^2$. Thus $K_S^2=K_S\cdot C=C^2=2$ and
then
$K_S\simnum C$.

If $p_g=1$, then $Z=0$. Let $D\in|K_S|$ be the unique effective divisor,
then
$2D=F_1+F_2$, where the $F_i$ are two fibres of $\phi_2$. If $F_1\ne F_2$,
then
the $F_i$ are multiple fibres and then $D\simnum 2F_0$, where $F_0$ is a
divisor. Which implies $D^2\ge 4$, a contradiction. If $F_1=F_2$, then
$D=F_1$ and thus $h^0(S, D)=2$, also a contradiction.

If $p_g=0$, because the $\pi_1$ of $S$ is a finite group  (Corary 5.8 of
\cite{1}), then
$h^1(S, K_S-C)=0$ by Lemma 4.8. Whereas we have $h^1(S, K_S-C)=h^1(S, C)=1$
by R-R, a
contradiction. Therefore we have $C^2=0$ and then $g(C)=2$.
\qed\enddemo

\proclaim{Proposition 4.10} Let $S$ be a minimal surface of type (II) or
(III).
Then $|2K_S|$ can not be composed of a rational pencil of genus two.
\endproclaim
\demo{Proof}
We refer to the proof of Proposition 3 and Theorem 3 of \cite{19}.
\qed\enddemo
\medskip

Thus we finally arrive at the following theorem of Xiao (Theorem 1 of \cite{19}). 

\proclaim{Theorem 4.11} Let $S$ be a projective surface of general type.
Then
$\fei{2}$ is generically finite if and only if $h^0(S,2K_S)>2$.
\endproclaim

\remark{\smc Acknowledgment}
I would like to thank the Abdus Salam International Centre for
Theoretical
Physics, Trieste, Italy for support during my visit.
I also wish to thank Prof. L. Ein,
Prof. M.S. Narasimhan and Prof. M. Reid for encouragement.
Thanks are also due to the hospitality of both Prof. C. Ciliberto and Prof.
E. Stagnaro who are very kind to invite me visiting their universities.
During the writing of this paper, a discussion with Prof. Ciliberto was quite helpful for me to organize the last section. 
This paper was finally revised while I was visiting the Uinversit$\ddot{a}$t G$\ddot{o}$ttingen, Germany as a post-doc fellow.
\endremark
\bigskip
\centerline{-----------------}
\centerline{\bf References}
\smallskip
\roster
\item"[1]" A. Beauville: {\it L'application canonique pour les surfaces de type
g\'en\'eral},  Invent. Math. {\bf 55}(1979), 121-140.
\item"[2]" X. Benveniste: {\it Sur les applications pluricanoniques des
vari\'et\'es
de type tr\'es g\'eg\'eral en dimension 3},  Amer. J. Math. {\bf
108}(1996), 433-449.
\item"[3]" E. Bombieri: {\it Canonical models of surfaces of general type},
 Publication I.H.E.S. {\bf 42}(1973), 171-219.
\item"[4]" M. Chen: {\it On pluricanonical maps for threefolds of general type},
J. Math. Soc. Japan {\bf 50}(1998), 615-621.
\item"[5]" -----: {\it A theorem on pluricanonical maps of nonsingular minimal threefold
of general type}, Chin. Ann. of Math. {\bf 19B}: 4(1998), 415-420.
\item"[6]" -----: {\it Complex varieties of general type whose canonical
systems are composed with pencils}, J. Math. Soc. Japan
{\bf 51}(1999), 331-335. 
\item"[7]" -----: {\it Kawamata-Viehweg vanishing and quint-canonical  map of a
complex threefold},  Comm. in Algebra {\bf 27}(1999), 5471-5486.
\item"[8]" O. Debarre: {\it Addendum, in\'egaliti\'es num\'eriques pour les
surfaces
de type g\'en\'eral},  Bull. Soc. Math. France {\bf 111}(1983), 301-302.
\item"[9]" M. Hanamura: {\it Stability of the pluricanonical maps of threefolds},
In:  Algebraic Geometry, Sendai, 1985 (Adv. Stud. in Pure Math. {\bf 10},
1, pp. 185-205).
\item"[10]" Y. Kawamata: {\it A generalization of Kodaira-Ramanujam's
vanishing theorem}, Math. Ann. {\bf 261}(1982), 43-46.
\item"[11]" J. Koll\'{a}r: {\it Higher direct images of dualizing
sheaves I}, Ann. of Math. {\bf 123}(1986), 11-42.
\item"[12]" S. Lee: {\it Remarks on the pluricanonical and adjoint linear series on projective threefolds}, Commun. Algebra {\bf 27}(1999), 4459-4476.
\item"[13]" K. Matsuki: {\it On pluricanonical maps for 3-folds of general
      type}, J. Math. Soc. Japan {\bf 38}(1986), 339-359.
\item"[14]" Y. Miyaoka: {\it The Chern classes and Kodaira dimension of
      a minimal variety},
     In: Algebraic Geometry, Sendai, 1985 (Adv. Stud. in Pure Math.
     {\bf 10}, 1987, pp.  449-476).
\item"[15]" M. Reid: {\it Canonical 3-folds}, in Journ\'ees de G\'eom\'etrie
Alg\'ebrique d'Angers (A. Beauville ed.), Sijthoff and Noordhof, Alphen ann
den Rijn, 1980, 273-310.
\item"[16]" I. Reider: {\it Vector bundles of rank 2 and linear systems
on algebraic surfaces}, Ann. of Math. {\bf 127}(1988), 309-316.
\item"[17]" E. Viehweg: {\it Vanishing theorems}, J. reine angew. Math.
{\bf 335}(1982), 1-8.
\item"[18]" G. Xiao: {\it L'irr\'egularit\'e des surfaces de type g\'en\'eral
dont
le syst\`eme canonique est compos\'e d'un pinceau}, Compositio Math.
{\bf 56}(1985), 251-257.
\item"[19]" -----: {\it Finitude de l'application bicanonique des surfaces de
type
g\'eg\'eral}, Bull. Soc. Math. France {\bf 113}(1985), 23-51.
\endroster
\enddocument